\documentclass[12pt]{amsart}

\newcommand{\pa}[2]{\ensuremath{\left(\begin{smallmatrix}{#1}\\{#2}
\end{smallmatrix}\right)}}
\newcommand{\bev}{\te{$\beta'$ even}}

\makeatletter
\@namedef{subjclassname@2020}{\textup{2020} Mathematics Subject Classification}

\usepackage{youngtab}

\makeatother
\newcommand{\run}{\te{run}}

\newcommand{\ssotknl}{\ssot_{kn}(\l)}

\newcommand{\ssyt}{\te{SSYT}}

\renewcommand{\lg}{*(lightgray) }
\newcommand{\bur}{\te{Bur}}
\newcommand{\pr}{\te{pr}}

\newcommand{\prs}{\te{pr}\,S}

\newcommand{\tbox}{\te{box}}

\newcommand{\mg}{\infty}

\newcommand{\qte}[1]{\q\te{#1}}

\renewcommand{\b}{\mathbf{b}}

\newcommand{\oll}{\ol{\el}}

\newcommand{\dsum}{\di\sum}

\newcommand{\supp}{\te{supp}\,}

\DeclareFontFamily{U}{mathx}{\hyphenchar\font45}
\DeclareFontShape{U}{mathx}{m}{n}{
      <5> <6> <7> <8> <9> <10>
      <10.95> <12> <14.4> <17.28> <20.74> <24.88>
      mathx10
      }{}
\DeclareSymbolFont{mathx}{U}{mathx}{m}{n}
\DeclareFontSubstitution{U}{mathx}{m}{n}
\DeclareMathAccent{\widecheck}{0}{mathx}{"71}

\newcommand{\fit}[1]{\mb{}\newline\resizebox{\textwidth}{!}{${#1}$}\newline}

\newcommand{\lr}{\longleftrightarrow}

\renewcommand{\kill}[1]{}
\newcommand{\dummy}[1]{\mbox{}}

\makeatletter
\newcommand{\xequal}[2][]{\ext@arrow 0055{\equalfill@}{#1}{#2}}
\def\equalfill@{\arrowfill@\Relbar\Relbar\Relbar}
\makeatother

\newcommand{\mto}{\mapsto}

\newcommand{\ku}{\ensuremath{\emptyset}}

\renewcommand{\k}{\ensuremath{\ol{\mathrm{P}}}}
\newcommand{\n}{\ensuremath{\bm{n}}}

\newcommand{\hou}[3]{{#1}\equiv {#2}\pmod{#3}}

\newcommand{\h}{\hline}

\newcommand{\kyo}[1]{{\eh{#1}}}

\newcommand{\ok}{}

\newcommand{\m}{\ensuremath{\infty}}

\renewcommand{\k}[1]{\ensuremath{\left({#1}\right)}}

\newcommand{\ds}{\dots}

\newcommand{\bca}{\begin{cases}}
\newcommand{\eca}{\end{cases}}

\newcommand{\des}{\textnormal{des}}

\renewcommand{\th}{\ensuremath{\theta}}

\newcommand{\del}{\ensuremath{\delta}}

\renewcommand{\ss}[3]{\ensuremath{\di\int_{#1}^{#2}{#3}\,dx}}

\newcommand{\bpic}{\begin{picture}}\newcommand{\epic}{\end{picture}}

\newcommand{\beda}{\begin{edaenumerate}}
\newcommand{\eeda}{\end{edaenumerate}}

\newcommand{\ret}[2]{\ensuremath{
\begin{pmatrix}{#1}\\  {#2}\end{pmatrix}
}}
\newcommand{\cd}{\cdots}

\newcommand{\sh}[1]{\shadowbox{#1}}
\newcommand{\st}{\strut}

\newcommand{\q}{\quad}

\newcommand{\too}{\longrightarrow}

\newcommand{\bq}{\begin{quote}}\newcommand{\eq}{\end{quote}}

\renewcommand{\sp}[1]{\ul{\ph{#1}}}
\renewcommand{\sp}[1]{\text{sp}}

\newcommand{\ti}{\times}
\newcommand{\ra}{\rangle}\newcommand{\la}{\langle}

\newcommand{\be}{\begin{enumerate}}\newcommand{\ee}{\end{enumerate}}
\newcommand{\bce}{\begin{center}}\newcommand{\ece}{\end{center}}
\newcommand{\bde}{\begin{description}}\newcommand{\ede}{\end{description}}
\newcommand{\bri}{\begin{flushright}}\newcommand{\eri}{\end{flushright}}
\newcommand{\bb}{\begin{block}}\newcommand{\eb}{\end{block}}
\newcommand{\bt}{\begin{thm}}\newcommand{\et}{\end{thm}}
\newcommand{\bpf}{\begin{proof}}\newcommand{\epf}{\end{proof}}
\newcommand{\bex}{\begin{ex}}\newcommand{\eex}{\end{ex}}
\newcommand{\bexr}{\begin{exr}}\newcommand{\eexr}{\end{exr}}
\newcommand{\bft}{\begin{fact}}\newcommand{\eft}{\end{fact}}
\newcommand{\brk}{\begin{rmk}}\newcommand{\erk}{\end{rmk}}
\newcommand{\ba}{\begin{align*}}\newcommand{\ea}{\end{align*}}
\newcommand{\bexe}{\begin{exe}}\newcommand{\eexe}{\end{exe}}

\newcommand{\bit}{\begin{itemize}}\newcommand{\eit}{\end{itemize}}

\newcommand{\bcm}{}
\newcommand{\ol}{\overline}\newcommand{\ul}{\underline}
\newcommand{\hf}{\hfill}

\newcommand{\cc}{\ensuremath{\mathbf{C}}}
\newcommand{\nn}{\ensuremath{\mathbf{N}}}

\newcommand{\zz}{\ensuremath{\mathbf{Z}}}

\newcommand{\bd}{\begin{defn}}\newcommand{\ed}{\end{defn}}
\newcommand{\bp}{\begin{prop}}\newcommand{\ep}{\end{prop}}

\newcommand{\eh}{\emph}\newcommand{\al}{\alpha}
\newcommand{\sub}{\subseteq}
\newcommand{\lam}{\lambda}
\newcommand{\fb}{\fbox}
\newcommand{\mb}{\mbox}
\newcommand{\te}{\text}\newcommand{\ph}{\phantom}
\newcommand{\wt}{\widetilde}\newcommand{\sm}{\setminus}
\newcommand{\ri}{\right}
\renewcommand{\l}{\left}

\newcommand{\di}{\displaystyle}\renewcommand{\a}{\ensuremath{\bm{a}}}
\renewcommand{\b}{\ensuremath{\bm{b}}}

\newcommand{\x}{\ensuremath{\bm{x}}}

\newcommand{\np}{\newpage}

\renewcommand{\b}{\beta}

\renewcommand{\a}{\alpha}

\renewcommand{\x}{\mathbf{x}}

\usepackage[dvipdfmx]{graphicx}
\usepackage[dvipsnames]{xcolor}
\usepackage{float,afterpage}
\usepackage{asymptote,layout}
\usepackage{wrapfig,epic}
\usepackage{amsmath,amssymb,cases,pict2e}

\usepackage{multicol}
\usepackage{tikz}
\tikzset{
    cell/.style={
        anchor=south west,
        draw,
        minimum size=1cm,
    },
}

\usepackage{geometry}
\usepackage{exscale,latexsym,bm}
\usepackage{amssymb,enumerate,amsmath,amsthm,amsfonts}
\usepackage{verbatim,fancybox,wasysym,fancyhdr,type1cm}
\usepackage[frame,all,poly,curve,knot,arrow]{xy}
\usepackage{ytableau}
\usepackage{colortbl}
\usepackage{boxedminipage}
\usepackage{multirow}

\graphicspath{{./figs/}}
\theoremstyle{definition}
\newtheorem{thm}{Theorem}[section]

\newtheorem{lem}[thm]{Lemma}

\newtheorem{prop}[thm]{Proposition}\newtheorem{cor}[thm]{Corollary}

\newtheorem{exr}[thm]{Exercise}
\newtheorem{ob}[thm]{Observation}

\newtheorem{ex}[thm]{Example}

\newtheorem{defn}[thm]{Definition}\newtheorem{rmk}[thm]{Remark}
\newtheorem{fact}[thm]{Fact}
\newtheorem{block}[thm]{}
\newtheorem*{exe}{Exercise}

\renewcommand{\a}{\alpha}

\renewcommand{\l}{\lam}
\renewcommand{\h}{\hline}
\renewcommand{\arraystretch}{1.5}

\newcommand{\innn}{\in\nn}

\usepackage{youngtab}

\newcommand{\ot}{\otimes}

\renewcommand{\l}{\lambda}
\renewcommand{\wt}{\text{wt}}
\renewcommand{\ep}{\varepsilon}
\newcommand{\yd}{\ydiagram}

\renewcommand{\m}{\mu}

\newcommand{\B}{\mathbf{B}}
\renewcommand{\del}{\text{del}}

\newcommand{\ssot}{\te{SSOT}}

\renewcommand{\st}{\text{st}\,}
\renewcommand{\ot}{\text{OT}}
\newcommand{\OT}{\text{OT}}

\renewcommand{\a}{\bm{a}}

\renewcommand{\sh}{\te{sh\,}}

\renewcommand{\ss}{\text{ss}}

\renewcommand{\a}{\al}
\renewcommand{\B}{\mathbf{B}}

\renewcommand{\a}{\alpha}
\newcommand{\wti}{\ensuremath{\text{wt}_{i}}}
\newcommand{\wts}{\ensuremath{\text{wt}^{*}\,}}
\renewcommand{\b}{\beta}

\renewcommand{\sh}{\te{sh\,}}

\renewcommand{\L}{\mathbf{L}}
\renewcommand{\B}{\mathbf{B}}
\renewcommand{\sh}{\te{sh}\,}

\renewcommand{\m}{\mu}

\renewcommand{\l}{\lambda}
\newcommand{\clmn}{c_{\l\m}^{\n}}

\renewcommand{\B}{\mathbf{B}}
\renewcommand{\bt}{\te{BT}}

\newcommand{\dest}{\te{des}\,T}

\newcommand{\parm}{\ensuremath{\te{Par}(m)}}
\newcommand{\parn}{\ensuremath{\te{Par}(n)}}
\renewcommand{\oll}{\ol{\lam}}
\newcommand{\olm}{\ol{\mu}}

\renewcommand{\del}{\partial}

\newcommand{\tref}{\te{ref}}

\newcommand{\cu}{\text{cu}}
\newcommand{\Cu}{\text{Cu}}

\renewcommand{\n}{\nu}
\newcommand{\Des}{\te{Des}}

\renewcommand{\wts}{\widetilde{S}}

\renewcommand{\L}{\Lambda}
\newcommand{\com}{\te{com}}
\newcommand{\Com}{\te{Com}}

\newcommand{\dst}{\te{dst}\,}
\newcommand{\step}{\te{step}\,}

\renewcommand{\flat}{\te{flat}}

\newcommand{\qsot}{\te{QYOT}}
\renewcommand{\del}{\te{del}}
\newcommand{\add}{\te{add}}
\renewcommand{\wti}{\widetilde}

\begin{document}

\title[generating function of SSOT]{Symmetry of the generating function of semistandard oscillating tableaux}


\author[M. Kobayashi]{Masato Kobayashi}
\author[T. Matsumura]{Tomoo Matsumura}
\thanks{corresponding author:\,Masato Kobayashi}
\author[S. Sugimoto]{Shogo Sugimoto}

\date{\today}                                       


\address{Masato Kobayashi\\
Department of Engineering\\
Kanagawa University, Rokkaku-bashi, Yokohama, Japan.}
\email{masato210@gmail.com}

\address{Tomoo Matsumura\\
Department of Natural Sciences\\
International Christian University, 
Tokyo, Japan.
}
\email{matsumura.tomoo@icu.ac.jp}

\address{Shogo Sugimoto\\
Department of Natural Sciences\\
International Christian University, 
Tokyo, Japan.
}
\email{shougo.sugimoto@gmail.com}

\maketitle
\begin{abstract}
H.Choi-D.Kim-S.J.Lee and S.J.Lee introduced a 
new kind of tableaux, semistandard oscillating tableaux (SSOT), around 2024 in the context of 
Lusztig $q$-weight multiplicities, KR crystals and King tableaux. 
In this paper, we study generating function of 
the SSOTs and its symmetry. 
First, we extend Gessel's and Assaf-Searles' expansion of a Schur function in terms of fundamental quasi-symmetric functions to our generating function. As a consequence, we show that it is $F$-positive. 
Further, we improve Sundaram's work on oscillating tableaux by proving that it is symmetric, Schur-positive, and has Saturated Newton polytope.
%
\end{abstract}


\tableofcontents

\ytableausetup{centertableaux}

%
%
%


\section{Introduction}

\subsection{Semistandard oscillating tableaux
and Gessel's theorem}

One of the important topics 
in algebraic/enumerative combinatorics 
is the theory of \kyo{Schur functions}. 
This idea is ubiquitous as they appear in the 
theory of the ring of symmetric functions,  Littlewood-Richardson rule, Jacobi-Trudi formula, to name a few. 
In particular, $s_{\l}$ is 
the generating function of 
semistandard Young tableaux of a partition shape $\l$. This is symmetric as Bender-Knuth proved in 1972 \cite{bek}.


Among other things, Gessel \cite{ge} in 1984 
introduced a \kyo{fundamental quasi-symmetric function}. He showed that, 
roughly speaking, a Schur function 
positively expands into a finite sum of 
those; 
Assaf-Searles in 2017 \cite{as} improved his theorem under polynomial setting. 
It is then natural to ask if 
analogous results hold for 
the generating function of other kind of tableaux.

Recently, H.Choi-D.Kim-S.J.Lee \cite{ckl} and S.J.Lee \cite{lee} 
introduced \kyo{semistandard} \kyo{oscillating} \kyo{tableaux} (SSOT) in context of Lusztig $q$-weight multiplicities, KR crystals and King tableaux. 
Inspired by them, the first and second authors \cite{km} found a bijective RSK correspondence 
with the set of pairs of King tableaux and SSOTs and introduced the generating function of SSOTs (which we call simply \kyo{SSOT function}) to prove a certain Cauchy identity. It turns out that 
we can extend the theorems 
by Gessel and Assaf-Searles \kyo{mutatis mutandis} 
for those functions. 
Furthermore, they are related to a variety of ideas  in this area such as descent, Burge array, Littlewood-Richardson coefficients, dominance order of partitions and even Newton polytopes.

\begin{rmk}
SSOT is the semistandard version of \kyo{oscillating tableaux} as 
Sundaram (1986) discussed with great details in  \cite{sun1}. 
Motivated by her work, Heo-Kwon \cite{hk}, 
Naito-Sagaki \cite{ns}, 
Kumar-Torres \cite{sat} and 
Rubey-Sagan-Westbury \cite{rsw} 
have studied the related representation theory of  symplectic groups. 
Although we stick to the combinatorics of tableaux, 
it should be possible to find connections 
to their work.
\end{rmk}


\subsection{Main results}

Our main results consist of a series of theorems as a typical type-A-to-other-type-extension where type A refers to the type of root system 
behind the scene.
Theorems \ref{t1} and \ref{t12} for OT/SSOTs 
generalize Gessel's and Assaf-Searles' theorems, respectively. 
In particular, an SSOT function is $F$-positive.
As a consequence, we prove as Theorems \ref{t31}, \ref{t41} that an SSOT function is 
symmetric and its homogeneous part is Schur- positive (as a finite sum of Schur functions).
Furthermore, Theorems \ref{t42}, \ref{t43}, \ref{t44} 
show partial orthogonality and 
linear independence, SNPness (Saturated Newton Polytope) of such functions. 
These are analogous properties 
which classical Schur functions satisfy.

\subsection{Organization}

Section 2 provides all fundamental ideas and definitions on semistandard Young tableaux, OT and SSOT. 
In Section 3, we revisit and improve Burge and Sundaram correspondences.
In Section 4, we show applications of 
symmetry of SSOT function 
to Saturated Newton Polytopeness with a simple combinatorial model.

\subsection{Acknowledgment.}

Authors would like to thank participants of FPSAC 2025 (Sapporo) 
and Hideya Watanabe for fruitful discussions.

\section{Calculus of type C tableaux}

We assume familiarity 
of basic combinatorics on SSYTs in Part I of Fulton's book \cite{fu}.
Letters $i, j, k, l, m, n$ are nonnegative integers 
unless otherwise specified; 
$\l, \m, \n$ are always partitions;
we use symbols $O, Q, R, S, T, U$ for tableaux.

\subsection{Gessel's quasi-symmetric function}

A \kyo{weak (resp. strong) composition} is a finite sequence 
of nonnegative (resp. positive) integers as we often write 
$a=(a_{1}, \ds, a_{m})$. 
For convenience, we identify compositions 
with ignoring 0's at the end as 
$(0, 2, 0)$ and $(0, 2)$. 

For simplicity, 
we sometimes abbreviate 
$a_{1}\cdots a_{m}$ to mean 
$(a_{1}, \ds, a_{m})$. 
The \kyo{size} of $a$ is 
\[
|a|=a_{1}+\cd+a_{m}.
\]
By $\Com$ (resp. $\ol{\Com}$) we mean the set of all 
weak (resp. strong) compositions. 
The \kyo{weight} (vector) of 
$a\in \Com$ is $\wt(a)=(\wt_{i}(a))$ 
where $\wt_{i}(a)$ is the number of $i$ in $a$.

Let $a, b\in \ol{\Com}$ and $c\in \Com$. 
Define $b=\flat(c)$ if we obtain $b$ with 
deleting all 0 parts of $c$. 
For example, 
if $c=(0, 2, 3, 0, 2)$, then $\flat(c)=(2, 3, 2)$. 

Say $b$ \kyo{refines} $a$ if 
 $a=(a_{1}, \ds, a_{m}),
b=(b_{1}, \ds, b_{n})$ and there exist some 
$i_{1}, \ds, i_{m}$ such that 
$1= i_{1}<\cd <i_{m}\le n$,
\[
a_{j}=b_{i_{j}}+\cd+b_{i_{j+1}-1}
\]
for each $j=1, \ds, m$ (and $i_{m+1}=n+1$). 
In particular, $a$ refines itself. Denote 
by $\tref(a)$ the set of all 
strong compositions which refine $a$.
For example, if $a=(2, 3)$, then 
\[
\tref{(a)}=
\{a, 221, 212, 2111, 113, 1121, 1112, 11111\}.
\]

Let 
$\x=(x_{1}, x_{2}, \ds)$ be infinitely many variables. 
To any composition $a$, 
associate the weight monomial 
$x^{a}=x_{1}^{a_{1}}x_{2}^{a_{2}}\cd$.


\begin{defn}
Let $a, b\in \ol{\Com}$. 
The \kyo{monomial quasi-symmetric function} for $b$ is 
\[
M_{b}(\x)=\sum_{
\substack{c\in \Com\\\flat(c)=b}
}
x^{c}.
\]
The \kyo{fundamental quasi-symmetric function} for $a$ is 
\[
F_{a}(\x)=\sum_{b\in \tref{(a)}}
M_{b}(\x).
\]
\end{defn}
To describe 
connections with such functions and SSYTs, we need another definition; a skew diagram is a \kyo{horizontal strip} if it contains at most one box at each column.

\begin{rmk}
In this article, a \kyo{box} always means \kyo{only a position} in a diagram. 
\end{rmk}
Let $B_{1}, B_{2}$ be boxes in a diagram. 
Write $B_{1}<_{nE} B_{2}$ if $B_{2}$ lies 
northEast to $B_{1}$, that is, 
if 
$B_{1}=(i_{1}, j_{1})$ and  
$B_{2}=(i_{2}, j_{2})$, 
then $i_{1}\ge i_{2}$ and $j_{1}<j_{2}$.

Suppose now 
$B=(B_{1}, \ds, B_{n})$ 
is a sequence of boxes and 
moreover there is an associated sequence 
$u=(u_{1}, \ds, u_{n})$ of postive integers 
such that $u_{j}$ is assigned to $B_{j}$.



\begin{defn}
Say $H=(B, u)$ as above is a \kyo{standard horizontal band} if 
\begin{itemize}
	\item $n>0$,
	\item $B_{1}<_{nE}\cd <_{nE}B_{n}$,
	\item $u_{j+1}-u_{j}=1$ for all $j$.
\end{itemize}
\end{defn}
By a \kyo{standard Young tableau} (SYT) $T$ of shape $\l$, we mean an SSYT such that 
each of $1, \ds, |\l|$ appears exactly once as an entry of $T$.
%

For such $T$, we can find 
a unique sequence of maximal standard horizontal bands 
$H_{1}, \ds, H_{k}$ 
such that 
they contain each entry of $T$ exactly 
once; 
if $u_{i}$ and $u_{j}$ appear in $H_{i}$ and $H_{j}$ $(i<j)$  respectively, then $u_{i}<u_{j}$; maximal means adding a box to $H_{i}$ 
cannot be any other standard horizontal band in $T$.
\begin{defn}
The \kyo{descent composition} of an SYT $T$ 
is the strong composition 
\[
\des\, T =\k{|H_{1}|, |H_{2}|, \ds, |H_{k}|}
\]
as $H_{i}$'s described just above. 
Define $\step T=k$.
%
\end{defn}

\begin{rmk}
In the context of quasi-symmetric functions, Young tableaux, permutations and words, 
this term \kyo{descent} often appears.
There are two essentially 
equivalent formulations, \kyo{descent composition} and \kyo{descent sequence (set)}. 
We will mainly deal with the former.
\end{rmk}

For example, if 
\ytableausetup{mathmode, boxsize=14pt}
$T=
\begin{ytableau}
	1&2&3&5\\
	4&7\\6
\end{ytableau}$,
then we find three maximal standard horizontal bands: 
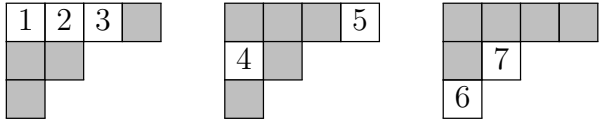
\[
\begin{ytableau}
	1&2&3&\lg{}\\
	\lg{}&\lg{}\\\lg{}
\end{ytableau}
\q
\q
\begin{ytableau}
	\lg{}&\lg{}&\lg{}&5\\
	4&\lg{}\\\lg{}
\end{ytableau}
\q
\q
\begin{ytableau}
	\lg{}&\lg{}&\lg{}&\lg{}\\
	\lg{}&7\\6
\end{ytableau}
\]
Thus, $\dest=(3, 2, 2)$ and $\step T=3$.


Let $\te{ST}(\l)$, $s_{\l}$ denote the set of 
all standard tableaux of shape $\l$ 
and the Schur function for $\l$.
\begin{thm}[Gessel {\cite{ge}}]\label{tg}
We have
\[
s_{\l}(\x)=
\sum_{T\in \te{ST}(\l)}F_{\dest}(\x).
\]
\end{thm}
Our first goal is to extend this theorem 
for OTs and SSOTs.





%
%
%
%



\subsection{OT and SSOT}


For partitions $\l$ and $\mu$, 
write $\l\lhd \mu$ if $\l\subset \mu$ and 
$|\mu|=|\l|+1$. 
For convenience, let 
\[
N(\l)=\{n\in\zz_{\ge0}\mid n\ge |\l|, 
\hou{n}{|\l|}{2}\}.
\]

\begin{defn}\ok 
An \emph{oscillating tableau} (OT or up-down tableau) $O$ of shape $\l$ of length $n$ $(\in N(\l))$ is a sequence of partitions 
\[
O=(O_{0}, O_{1}, \ds, O_{n})
\]
such that $O_{0}=\ku$, $O_{n}=\l$ and for each $j$, 
$O_{j}\lhd O_{j+1}$ or $O_{j+1}\lhd O_{j}$. 
Denote by $\OT_{}(\l)$ the set of OTs of shape $\l$.
\end{defn}
For example, 
\ytableausetup{mathmode, boxsize=12pt}\[
O=\k{
\ku, 
\yd{1}, \yd{1,1}, \yd{2,1}, \yd{3, 1}, 
\yd{2, 1}, \yd{2}
}
\]
is an OT of length 6.

It is helpful to identify an OT 
with a single set-valued tableau by recording 
all history on addition/deletion of boxes. 
For example, we identify $O$ just above with  
$\ytableausetup{mathmode, boxsize=15pt}
\begin{ytableau}
	1&3&45\\
	26\\
\end{ytableau}.$
Define the \kyo{profile} $\te{pr}\,O=(1, \ds, n)$ 
for any OT of length $n$. 
Let $B_{j}$ be the box 
indicating the difference between $O_{j}$ and $O_{j+1}$. 
The \kyo{box-trace} of $O$ is 
$\te{box}\,O=(B_{j})_{j=1}^{n}$.
We may further extract crucial information from an OT with ``assembling" certain consecutive partitions. 

\begin{defn}
Let 
$O=(O_{j})$ be an OT of length $n$.
Define the \kyo{descent set (sequence)} 
$\Des\,O$ 
and \kyo{descent composition} 
$\des\,O$ of $O$ as follows. 
First, let $d(O)$ be the sequence 
\[
d(O)=
(1\le d_{1}\le 
d_{2'}\le d_{2}\le d_{3'}\le \cd \le d_{k-1}\le n-1) 
\]
with $d_{0}=d_{1'}=0, 
d_{k'}=d_{k}=n$ and $d_{i}, d_{{i}'}$ are numbers satisfying the following.

\begin{itemize}
\item 
$O_{d_{{i}'}}\lhd
O_{d_{{i}'}+1}\lhd
\cd \lhd O_{d_{i}}$
and $O_{d_{i}}/O_{d_{i'}}$
is a horizontal strip such that 
\[
O_{j+1}/O_{j}<_{nE}O_{j+2}/O_{j+1}
\q \te{for all $j$, $d_{i'}\le j\le d_{i}-2$}
\]
and moreover one of the following holds.
\begin{enumerate}
\item[(i)] 
$O_{d_{i}}\not\!\!\lhd\,\, O_{d_{i}+1}$.
\item[(ii)] 
$O_{d_{i}}\lhd O_{d_{i}+1}$ 
and 
\[
O_{d_{i}}/O_{d_{i}-1}\not<_{nE}O_{d_{i}+1}
/O_{d_{i}}.
\]
\end{enumerate}
\item 
$O_{d_{i}}\rhd
O_{d_{i}+1}\rhd
\cd \rhd O_{d_{{(i+1)'}}}$ 
and 
$O_{d_{i}}/O_{d_{{(i+1)'}}}$ is a horizontal strip such that 
\[
O_{j+1}/O_{j+2}<_{nE}O_{j}/O_{j+1}
\q \te{for all $j$, $d_{i}\le j\le d_{(i+1)'}-2$}
\]
and moreover one of the following holds.
\begin{enumerate}
\item[(i)] 
$O_{d_{(i+1)'}}\not\!\!\rhd\,\, O_{d_{(i+1)'}+1}$.
\item[(ii)] 
$O_{d_{(i+1)'}}\rhd\,\, O_{d_{(i+1)'}+1}$
and 
\[
O_{d_{(i+1)'}}/O_{d_{(i+1)'}+1}\not<_{nE}
O_{d_{(i+1)'}-1}
/O_{d_{(i+1)'}}.
\]
\end{enumerate}

\end{itemize}
Now define 
\[
\Des\,O=\{d_{i}\mid {1\le i\le k-1}\}, \q 
\des\,O=(d_{i}-d_{i-1})_{1\le i\le k}.
\]
Further, set $\step\,O=k$.
\end{defn}

Define 
$\ot_{kn}(\l)$ 
to be the subset of 
$\ot(\l)$ consisting of OTs of step $k$ and length $n$.
%

\begin{ex}
Let 
\ytableausetup{mathmode, boxsize=8pt}
\[
O=
\k{
\begin{array}{ccc|c|cccc|c}
	O_{0}&O_{1}&O_{2}&O_{3}&O_{4}&O_{5}&O_{6}&O_{7}&O_{8}\\
\ku&\yd{1}&\yd{2}&\yd{2,1}&\yd{1,1}&\yd{1}&\yd{1,1}
&\yd{2,1}&\yd{2,2}\\
\end{array}
}.
\]
Observe that $n=8$, $k=\step (O)=4$, 
The vertical lines indicate 
$d_{1}=2, d_{2}=3, d_{3}=7$. 
Also, $d_{2'}=2, d_{3'}=5, d_{4'}=7$.
Therefore, 
\[
\Des\, O=\{2, 3, 7\}, 
\des\, O=(2, 1, 4, 1).
\]
\end{ex}

\begin{rmk}
If we think 
$\ssyt(\l)\sub 
\ot(\l)$, then our descent set and composition are compatible with the usual ones for SSYTs. 
\end{rmk}

\begin{defn}
Let $n\in N(\l), k\in\zz_{\ge0}$.
A \kyo{semistandard oscillating tableau} (SSOT) of shape $\l$ of length $n$ with $k$ steps 
is a sequence of partitions 
$S=(S^{1}, S^{'2}, S^{2}, S^{'3}, \ds)$ which satisfies the following.
\begin{enumerate}
\item 
$S^{k}=S^{' k+1}=S^{k+1}=\cd=\l.$
\item $S^{i}\supseteq S^{'i+1}$ $(1\le i\le k-1)$
and $S^{'i}\sub S^{i}$ ($2\le i\le k$). Moreover, 
each of $S^{i} \backslash S^{'i+1}$ and 
$S^{i} \backslash S^{'i}$ 
is a (possibly empty) horizontal strip.
\item We have 
\[
n=
\sum_{i=1}^{k-1}
|S^{i} \backslash S^{'i+1}|+
\sum_{i=2}^{k}
|S^{i} \backslash S^{'i}|.
\]
\end{enumerate}
We write $\step(S)=k$.
Denote by $\ssot_{kn}(\l)$ the set of all such SSOTs. 
\end{defn}

By definition, such a sequence becomes stable after a finite step. 
Thus, we may ignore its tail part.
Unless otherwise specified, 
the letter $S$ means an SSOT in what follows.

\begin{rmk}
This definition is slightly different from Lee  \cite{lee}. An SSOT in his style has a step 
with adding boxes first and deleting them next.
The purpose of this technical modification is to 
construct the bijection for type C RSK correspondence by \kyo{Berele insertion} and a Knuth array \cite{km}; use of letter $u$ comes from the top row in such arrays.
\end{rmk}

\begin{ex}\label{ex1}
\[
\ytableausetup{mathmode, boxsize=15pt}
S=
\k{
\yd{2}, 
\yd{1}, 
\yd{3,1}, 
\yd{1,1},
\yd{2,1}
}
\]
is an SSOT.
\end{ex}

There is a convenient way 
to compactly express an SSOT by a single multi-set-valued tableau. 
We think that we add/delete one box at each  \kyo{sub}step as an $OT$.
Every time a box is added or deleted, we record its  \kyo{step number} into a tableau (of possibly larger than the final shape) at the same position. 
For example, we can encode 
the above $S$ as 
\[
\ytableausetup{mathmode, boxsize=15pt}
\k
{\begin{ytableau}
	1&1
\end{ytableau}, 
\begin{ytableau}
	\mb{}&\lg{2}
\end{ytableau}, 
\begin{ytableau}
	\mb{}&2&2\\2
\end{ytableau}, 
\begin{ytableau}
	\mb{}&\lg{3}&\lg{3}\\
	\mb{}
\end{ytableau}, 
\begin{ytableau}
	\mb{}&3\\
	\mb{}
\end{ytableau}}
\]
where the shaded boxes indicate deleted ones from the previous diagram. 
All together, we identify $S$ with 
the multi-set-valued tableau 
\[
\ytableausetup{mathmode, boxsize=38pt}
\begin{ytableau}
	1&12233&{23}\\
	2\\
\end{ytableau}
\]
as we just write down elements of a multi-set in each box.
Let $\|S\|$ mean the \kyo{content} of $S$, 
the multiset of all letters in $S$, say $\{u_{1}, \ds, u_{n}\}$.
It is essential to understand that there is always the canonical total order on $\|S\|$. 
Define a total order on $\|S\|$ 
by assigning positive numbers to each entry 
in $S^{i}$ from \kyo{left}, and one in $S^{'i}$ from \kyo{right} for $i=1, 2, \ds, \step S$. 
For example, 
\[
\|S\|=\{1_{(1)}, 1_{(2)}, 2_{(3)}, 2_{(4)}, 2_{(5)}, 2_{(6)}, 3_{(7)}, 3_{(8)}, 3_{(9)}\}
\]
with $S$ just as above.

\begin{defn}
Define the \eh{profile} of $S$ ($\te{pr\,}S$) of length $n$ 
to be the weakly increasing 
sequence consisting 
of all entries of $\|S\|$ with the total order.
If
\[
\te{pr\,}S=(u_{1}, \ds, u_{n}),
\]
then, the \kyo{box-trace} of $S$ is a sequence of boxes 
\[
\te{box}\, S=(B_{1}, \ds, B_{n})
\]
such that $B_{j}$ is the added/deleted box corresponding to $u_{j}$.
\end{defn}

\begin{ob}\label{ob1}
For SSOTs $R$ and $S$, we have:
$R=S$ if and only if 
$\te{pr\,}R=
\te{pr\,}S$ and 
$\te{box\,}R=
\te{box\,}S$.
\end{ob}

 
%

%
%
%

\subsection{Extension of Gessel's theorem}

Toward a proof of an extension of 
Gessel's theorem, 
let us deal with OTs and SSOTs 
of all steps together.

Let 
\[
\ot_{n}(\l)=
\bigsqcup_{k\ge0}\ot_{kn}(\l), \q 
\ssot_{n}(\l)=
\bigsqcup_{k\ge0}\ssotknl.
\]

\begin{defn}
Let $S\in \ssot_{n}(\l)$ 
with 
$\pr\, S=(u_{1}, \ds, u_{n})$.
Define $\st S$, its \kyo{standardization}, 
as the tableau obtained by replacing all $u_{j}$ by $j$.
This defines a map
\[
\st:\ssot_{n}(\l)\to \ot_{n}(\l).
\]
%
(See Example \ref{ex1} below).
\end{defn}

Say SSOTs $R$ and $S$ are 
\kyo{standard equivalent} 
if $\st R=\st S$. 
This is an equivalent relation on $\ssot_{n}(\l)$. 
It then makes sense to talk about its equivalent class $\ssot(O):=\st^{-1}(O)$ so that
\[
\ssot_{n}(\l)=\bigsqcup_{O\in \ot_{n}(\l)}\ssot(O).
\]


\begin{defn}
Define $\Des(S)=\Des(\st S)$ and 
$\des(S)=\des(\st S)$. 
\end{defn}

The standardization 
preserves des, length, box-trace, 
shape while it may reduce the step number:
$\step(\st S)\le \step (S).$


\begin{defn}
A \kyo{run} $u=(u_{j})$ is a finite 
sequence of positive integers and 
``$|$" with all of the following.
\begin{enumerate}
\item There is at most 
one $|$ between those adjacent integers.
\item It begins and ends at an integer.
\item $u_{j}\le u_{j+1}$ for all $j$.
\item If $u_{j}|u_{j+1}$, 
then $u_{j}<u_{j+1}$.
\end{enumerate}
Say $j$ is a \kyo{descent} of $u$ if 
$u_{j}|u_{j+1}$. 
The \kyo{weight} of $u$ is 
the composition $\wt\, u=(\wt_{i}(u))$ 
where $\wt_{i}(u)$ is the number of $i$'s in $u$ for each $i\le \max (u_{j})$. 
The length of a run $u$ is the number of integers 
and $\step(u)=$ the number of $|+1$.
\end{defn}

For example, 
$u=1112|3|55$ is a run 
of length 7 and of step 3 with $\wt\, u=(3, 1, 1, 0, 2)$. 
Notice that $u_{j}<u_{j+1}$ is possible 
even if there is no $|$ between 
$u_{j}$ and $u_{j+1}$.

\begin{defn}
Let $S$ be an SSOT such that $\pr\, S=(u_{j})$. 
The total run of $S$ $(\run\, S)$ is 
the sequence with the 
symbol $|$ between $u_{j}$ and $u_{j+1}$ inserted   
whenever $j\in \Des\,S$.
Similarly, define $\run\, O$ for an OT $O$.
\end{defn}


\begin{ex}
\label{ex1}
\ytableausetup{mathmode, boxsize=10pt}
Let 
$\l=
\yd{4, 1}$,
\ytableausetup{mathmode, boxsize=14pt}
\[
O=
\k{
\begin{array}{ccc|cccccc}
	O_{1}&O_{2}&O_{3}&O_{4}&O_{5}&O_{6}&O_{7}&O_{8}&O_{9}\\
\yd{1}&\yd{2}&\yd{3}&\yd{2}&\yd{1}&\yd{1,1}
&\yd{2,1}&\yd{3,1}&\yd{4,1}      \\
\end{array}
},
\]
\[
S=
\k{
\begin{array}{c|cc|cc}
	S^{1}&S^{'{2}}&S^{2}&S^{'{3}}&S^{3}      \\
\yd{3}	&\yd{1}&\yd{2,1}&\yd{2,1}&\yd{4,1}      \\
\end{array}
}.
\]
As set-valued tableaux, 
\ytableausetup{mathmode, boxsize=25pt}
\[
O=
\begin{ytableau}
	1&257&348&9\\
	6
\end{ytableau}
,
\q 
S=
\begin{ytableau}
	1&122&123&3\\
	2
\end{ytableau}.
\]
Observe that $\st S=O$, 
$\step\,O=2$ and $\step S=3$.
\[
\run\, O=123|456789, \q
\run\, S=111|222233,
\]
\[
d(O)=(d_{1}, d_{2'}, d_{2})=(3, 5, 9), 
\Des\, O=\{3\}, 
\des\, O=(3, 6).
\]
For understanding the proof of Theorem \ref{t1}, it is helpful to see some other examples 
of SSOTs $R$ as follows.
\[
\begin{array}{lclccccc}
R& &\run R\\\h
(O_{3}, O_{5}, O_{9})&&111|222222\\
(O_{3}, O_{4}, O_{4}, O_{5}, O_{9})&&111|233333\\
(O_{1}, O_{1}, O_{1}, O_{1}, O_{3}, O_{5},  O_{9})&&133|444444\\
\end{array}
\]
\end{ex}

\begin{defn}
For $O\in \ot_{n}(\l)$, define 
$\com:\ssot(O)\to \Com$ by 
$\com(S)=\wt (\pr S)$.
\end{defn}

\begin{lem}\label{linj}
$\com$ is injective and weight-preserving.
\end{lem}
\begin{proof}
Suppose 
$\com(R)=\com(S)$ 
for $R, S\in \ssot(O)$.
Then 
$\wt (\pr R)=\wt (\pr S)$ 
implies $\pr R=\pr S$ and 
$\tbox\,  R=\tbox\,O=\tbox\, S$.
Thus, $R=S$ as in Observation \ref{ob1}. 
Weight-preservingness is clear.
\end{proof}

Now let us clarify the basic idea on a set of compositions and the flat operation. 
Consider $C\sub \Com$. 
For $c, d\in C$, the relation 
$c\sim d\iff \flat(c)=\flat(d)$ 
defines an equivalent relation on $C$. Say 
\[
[c]=\{d\in C\mid d\sim c\}
\]
denotes an equivalent class 
and put $\ol{C}=C\cap \ol{\Com}$.
 Since $b=\flat(b)$ for $b\in \ol{C}$, it is natural to take such compositions as representatives:
\[
C=\bigsqcup_{b\in \ol{C}}\,[b].
\]

\begin{lem}\label{lcom}
For $O\in \ot_{n}(\l)$, 
$k'=\step O$, 
let $a_{O}=\des\, O
=(a_{1}, \ds, a_{k'})$ 
and $C=\com(\ssot(O))$.
Then 
$\ol{C}=\tref(a_{O})$.
\end{lem}

\begin{proof}
Suppose 
$b\in \ol{C}=
\ol{\com(\ssot(O))}$. 
Then 
$b=\com(S) $
for some $S\in \ssot(O)$ and 
 $b$ is strong.
Since $\com$ is injective as just shown above 
so that 
$S$ is unique, say $k=\step S$ and 
write 
\[
b=(b_{1}, \ds, b_{k}), b_{i}>0.
\]
It follows from $S\in \ssot(O)$ 
that 
\[
\des\, S=\des\, O=a_{O}, k'=\step O\le \step S=k
\]
and 
\[
\run\, S=\underbrace{1^{b_{1}}\cd}_{
\te{$a_{1}$ letters}} |\cd|\underbrace{\cd k^{b_{k}}}_{\te{$a_{k'}$ letters}}.
\]
This means exactly $b\in \tref(a_{O})$.

To prove the converse, 
suppose $b\in \tref(a_{O})$. 
Again, write 
\[
b=(b_{1}, \ds, b_{k}), b_{i}>0, k\ge k'.
\]
Now choose 
$S\in \ssot(O)$ such that 
\[
\pr S=(1^{b_{1}}, \ds, k^{b_{k}})
\]
which is possible 
since $b$ refines $a_{O}
=(a_{1}, \ds, a_{k'})$, that is, again, 
\[
\run\, S=\underbrace{1^{b_{1}}\cd}_{
\te{$a_{1}$ letters}} |\cd|\underbrace{\cd k^{b_{k}}}_{\te{$a_{k'}$ letters}}.
\]
With this choice of $S$, we have 
\[
\com S=\wt (\pr S)=\wt(1^{b_{1}}, \ds, k^{b_{k}})=b
\]
so that 
$b\in \ol{\com(\ssot(O))}$. 
%
\end{proof}

For $S\in \ssot_{n}(\l)$, 
let 
$\x^{S}=\prod_{i\in \|S\|}x_{i}$ be its weight monomial.
\begin{defn}
Let $n\in N(\l)$. Define the \kyo{SSOT function} of degree $n$ associated to $\l$ by 
\[
\ss_{\l, n}(\x)=\sum_{S\in \ssot_{n}(\l)}\x^{S}.
\]
For convenience, 
let $\ss_{\l, n}(\x)=0$ if $n\not\in N(\l)$.
\end{defn}
In particular, $\ss_{\l, |\l|}(\x)=s_{\l}(\x)$.

We now generalize Gessel's theorem.
\begin{thm}\label{t1}
Let $\x=(x_{1}, \ds, )$ be infinite variables. 
Then we have
\[
\ss_{\l, n}(\x)=
\sum_{O\in\ot_{n}(\l)}F_{\des\, O}(\x).
\]
\end{thm}

\begin{proof}
For convenience, let 
$C=\com(\ssot(O))$ for $O\in \ot_{n}(\l)$.
We have \[
\ss_{\l, n}(\x)=
\sum_{S}\x^{S}
=
\sum_{O}\sum_{S\in\ssot(O)}
\x^{S}
=
\sum_{O}
\sum_{c\in C}
\x^{c}
\]
\[
=
\sum_{O}
\sum_{b\in \tref{(a_{O})}}
\sum_{\flat(c)=b}\x^{c}
=
\sum_{O}
\sum_{b\in \tref{(a_{O})}}
M_{b}
=
\sum_{O\in\ot_{n}(\l)}
F_{\des\, O}(\x)
\]
where the third equality 
follows from 
Lemma \ref{linj}, 
the fourth from 
Lemma \ref{lcom}.

\end{proof}
\begin{cor}
$\ss_{\l, n}(\x)$ is $F$-positive.
\end{cor}

\subsection{Extension of Assaf-Searles' theorem}

Now consider the generating function of 
SSOTs of steps at most $k$, i.e., let 
$x_{k+1}=x_{k+2}=\cd=0$. 
In this case, infinitely many zero terms 
(for $O\in \ot_{n}(\l), \step O>k$) 
appear in the summation of Theorem \ref{t1}.
Following Assaf-Searles \cite{as}, 
we improve this situation a bit 
by choosing nice representatives 
of standard equivalent classes from 
$\ssot(\l)$ itself, rather than $\ot(\l)$.
%
%


Say an SSOT $Q$ is \kyo{quasi-Yamanouchi} if 
$\wt\, Q=\des\, Q$.
By $\qsot_{n}(\l)$ we mean 
the set of such tableaux 
of length $n$ and shape $\l$. 
Note that each 
$\ssot^{-1}(O)$ contains exactly one 
quasi-Yamanouchi tableau.

\begin{defn}[destandadization]
Define 
\[
\dst:\ssot_{n}(\l)\to \qsot_{n}(\l), \q 
\dst S=Q
\]
if $Q$ is the unique quasi-Yamanouchi tableau 
in $\ssot^{-1}(\st S)$.
\end{defn}


\begin{ex}
\[
O=
\begin{ytableau}
	1&2&347\\
	5&6
\end{ytableau}, \q
Q=
\begin{ytableau}
	1&1&122\\
	2&2
\end{ytableau}, \q
S=
\begin{ytableau}
	1&1&236\\
	3&4
\end{ytableau}
\]
satisfy $\st S=\st Q=O$ and 
$\dst S=Q$.
\end{ex}

\begin{ob}\label{lpre}
Maps $\st$ and $\dst$ restrict to des- and step-preserving bijections 
$\ot_{n}(\l)\cong\qsot_{n}(\l)$.
\end{ob}
When we consider 
only SSOTs of step at most $k$, 
we get a consequence of Theorem \ref{t1}.
Let $\ot_{\le k, n}(\l)$ and $
\qsot_{\le k, n}(\l)$ be the sets of 
corresponding tableaux of step at most $k$.
\begin{thm}\label{t12}
Let $\x_{k}=(x_{1}, \ds, x_{k})$ be finite variables. 
Then 
\[
\ss_{\l, n}(\x_{k})=
\sum_{Q\in \qsot_{\le k, n}(\l)}F_{\des \,Q}(\x_{k}).
\]
Moreover, each $F_{\des \,Q}(\x_{k})$ is nonzero.
\end{thm}
\begin{proof}
The equality follows from 
Theorem \ref{t1} and Observation \ref{lpre}.

To prove the second assertion, 
let $Q\in \qsot_{\le k, n}(\l)$. 
Then, 
\[
\des\, Q=a=(a_{1}, \ds, a_{j}),
\qte{and}\q
\run\, Q=1^{a_{1}}|\cd |j^{^{a_{j}}}
\]
for some $a\in \ol{\Com}$ and $j$ $(\le k)$ so that 
\[
F_{\des \,Q}(x_{1}, \ds, x_{k})
=
\sum_{b\in \tref{(a)}}M_{b}
=
x_{1}^{a_{1}}\cd x_{j}^{a_{j}}+\cd
\]
is nonzero.
\end{proof}

\begin{ex}
Table \ref{ta1} shows all 14 quasi-Yamanouchi SSOTs 
on $\l$ of length 5 of step at most 3 and their  runs; indeed, we exhausted all cases considering a  unique box containing three letters.
{
\renewcommand{\arraystretch}{3.5}
\begin{table}[h!]
\ytableausetup{mathmode, boxsize=8pt}
\caption{$\te{QYOT}_{\le 3, 5}\k{
\yd{2,1}}$}
\label{ta1}
\ytableausetup{mathmode, boxsize=28pt}
\begin{center}
\begin{tabular}{llllllccc}\h
\rule[-10mm]{0mm}{0mm}
\scalebox{0.82}
{$\begin{ytableau}
	1&1&12\\2\\
\end{ytableau}$}
	&$111|{2}2$	&	
&
\scalebox{0.82}
{$\begin{ytableau}
	1&122\\2\\
\end{ytableau}$}
	&$11|{2}22$	&	\\\h	
\rule[-10mm]{0mm}{0mm}
\scalebox{0.82}{$\begin{ytableau}
	1&1&13\\2\\
\end{ytableau}$}
	&$111|2|{3}$	&	&
\scalebox{0.82}{$\begin{ytableau}
	1&1&23\\2\\
\end{ytableau}$}
	&$11|22|{3}$	&	\\\h
\rule[-10mm]{0mm}{0mm}
\scalebox{0.82}{$\begin{ytableau}
	1&1\\2&23\\
\end{ytableau}$}
	&$11|22|{3}$	&	&
\scalebox{0.82}{$\begin{ytableau}
	1&122\\3\\
\end{ytableau}$}
	&$11|{2}2|3$	&	\\\h
\rule[-10mm]{0mm}{0mm}
\scalebox{0.82}{$\begin{ytableau}
	1&1\\233\\
\end{ytableau}$}
	&$11|2|{3}3$	&	&
\scalebox{0.82}{$\begin{ytableau}
	1&133\\2\\
\end{ytableau}$}
	&$11|2|{3}3$	&	\\\h
\rule[-10mm]{0mm}{0mm}
\scalebox{0.82}{$\begin{ytableau}
	1&2&23\\2\\
\end{ytableau}$}
	&$1|222|{3}$	&	&
\scalebox{0.82}{$\begin{ytableau}
	122&2\\3\\
\end{ytableau}$}
	&$1|{2}22|3$	&	\\\h
\rule[-10mm]{0mm}{0mm}
\scalebox{0.82}{$\begin{ytableau}
	1&2\\233\\
\end{ytableau}$}
	&$1|22|{3}3$	&	&
\scalebox{0.82}{$
\begin{ytableau}
	1&233\\2\\
\end{ytableau}$}
	&$1|22|{3}3$	&	\\\h
\rule[-10mm]{0mm}{0mm}
\scalebox{0.82}{$\begin{ytableau}
	122&3\\3\\
\end{ytableau}$}
	&$1|{2}2|33$	&	&
\scalebox{0.82}{$
\begin{ytableau}
	1&3\\233\\
\end{ytableau}$}
	&$1|2|{3}33$	&	\\\h
\end{tabular}
\end{center}
\end{table}}

From this, we 
find an $F$-expansion
\[
\ss_{\l, 5}(\x_{3})=
F_{32}+F_{23}+F_{311}+3F_{221}+2F_{212}+
2F_{131}+3F_{122}+F_{113}
\]
without terms such as $F_{2111}$.
\end{ex}

\np

%
%
%
%
%
%
%
%
%


\np

\section{Burge, Sundaram correspondences}

In this section, we concern the 
\kyo{full SSOT function} 
\[
\ss_{\l}(\x)=\sum_{n}\ss_{\l, n}(\x)
\]
with infinite variables.

\begin{rmk}
Our ideas in this section 
are based on Sundaram's work \cite{sun1} on King tableaux. 
Moreover, Schumann-Torres' discussion  \cite[Sections 4, 5]{sct} 
with Kashiwara-Nakashima tableaux 
is close to ours. 
The point here is that we can ``forget" 
everything about those tableaux and 
may just focus on SSOTs.
\end{rmk}

\subsection{Burge correspondence}

First, we review Burge correspondence \cite{bu}. A partition is \kyo{even} if each of its rows contains even number of boxes.


\begin{defn}
An array 
$L=
\left(\begin{array}{ccccc}
	j_{1}&\cd&j_{r}\\
	i_{1}&\cd &i_{r}
\end{array}\right)$ ($i_{l}, j_{l}\in \nn$)
is \kyo{lexicographic} \kyo{(Knuth)} 
if 
$j_{l}\le j_{l+1}$ and 
whenever $j_{l}=j_{l+1}$, then 
$i_{l}\le i_{l+1}$ for all $l\in [r]$.
It is \kyo{Burge} if moreover 
$j_{l}> i_{l}$ for all $l$.
An SSYT of shape $\b$ is a \kyo{Burge tableau} if 
$\b'$, the conjugate of $\b$, is even. 
Let $\B(2r)$ denote the set of Burge arrays with $r$ columns and $\bt(2r)$ the set of all 
SSYTs of Burge shape with $2r$ boxes.
\end{defn}

For convenience, we provide more notation and symbols.
For an array 
$L$ as above and $\pa{j}{i}$, 
$L\oplus \pa{j}{i}$ means the lexicographic array 
with rearranging all pairs of $L$ and 
$\pa{j}{i}$. Also, $2L$ means the \kyo{symmetrization} of $L$, 
the lexicographic array with all pairs 
$\pa{j_{l}}{i_{l}}$ in $L$ and $\pa{i_{l}}{j_{l}}$ 
with the same multiplicity.

\begin{fact}[Burge correspondence]
There is a bijection 
$\bur:\B(2r)\cong \bt(2r)$.
To be more precise, 
if 
$L=
\left(\begin{array}{ccccc}
	j_{1}&\cd&j_{r}\\
	i_{1}&\cd &i_{r}
\end{array}\right)$ 
is Burge 
and 
$2L=
\left(\begin{array}{ccccc}
	u_{{1}}&\cd&u_{{2r}}\\
	v_{{1}}&\cd &v_{{2r}}
\end{array}\right)$ 
is its symmetrization, then 
\[
\bur(L)=P(v_{{1}}\cd v_{{2r}})
\in \bt(2r)
\]
where $P(w)$ is the inserting tableau 
for a word $w$ with RSK correspondence (See 
\cite{fu} and \cite[p.38, Lemma 3.33]{sun1}).
\end{fact}

Often, we identify such a Burge array and Burge tableau under this correspondence.

\begin{ex}
\ytableausetup{mathmode, boxsize=15pt}
If 
$L=\left(\begin{array}{ccccc}
	4&4&7\\
	2&3&2
\end{array}\right)$, 
then 
$2L=\left(\begin{array}{cccccc}
	2&2&3&4&4&7\\
	4&7&4&2&3&2\\
\end{array}\right)$, 
\[
\bur
L=
P(474232)=
\begin{ytableau}
	2&2\\
	3&4\\
	4\\
	7\\
\end{ytableau}.
\]
\end{ex}


\subsection{Sundaram correspondence}

Sundaram discussed a certain 
correspondence 
between words and triples 
of a King tableau, an OT and 
a Burge array. 
As an improvement of her idea, 
we will construct a bijective correspondence 
$\th(S)$ as a pair $(L(S), T(S))$ 
from an SSOT $S$ alone. 
For this purpose, let us recall 
two operations, 
\kyo{reverse column-inserting} and 
\kyo{column-unbumping} for a tableau 
from 
the context of SSYTs.
Let 
$c_{j}, x\innn$ and 
$C=
\begin{ytableau}
	c_{1}\\
	\vdots\\
	c_{r}
\end{ytableau}$ be a column.
\begin{defn}
The reverse column-insert 
$x$ into $C$ is to 
find maximal $j$ such that $c_{j}\le x$ and 
then replace $c_{j}$ by $x$ 
and bump out $c_{j}$. 
\end{defn}

\begin{defn}
The column-unbumping of 
$c_{r}$ (the last entry of a column) 
from $C$ is an operation to bump out 
$c_{r}$.
\end{defn}
More generally, 
we can column-unbump 
an entry $x$ in a box $B$ at a corner of a tableau $T$ with the following procedure.
First, column-unbump $x$ from 
the column containing it. 
Then we get the bumped out number, say $x_{1}$ $(\le x)$. 
Then reverse column-insert $x_{1}$ to the previous column. Continue this algorithm 
back to the first one.
Finally, we obtain the SSYT with one less box (at $B$)  
and one bumped out number.
%

\begin{defn}
Let $T, B, x$ as above. 
By $\Cu(T, B)$ we mean the resulting tableau 
after \kyo{column-unbumping} $x$ at $B$ from $T$. 
Also, by $\cu(T, B)$ we mean the number bumped out from the first column by this operation.
\end{defn}
Whenever we want to emphasize 
an entry $x$ in $B$, we write 
$\Cu(T, B, x)$ and $\cu(T, B, x)$.
Note that $\cu(T, B, x)\le x$.
Also, if $B$ is at an outside corner of a tableau $T$, 
$T\oplus (B, x)$ means 
the tableaux with $T$ and $x$ at $B$.
Therefore, 
\[
T_{2}={T_{1}}\oplus (B, x)=i\to T_{1}
\iff T_{1}=\Cu(T_{2}, B, x), i=\cu(T_{2}, B, x).
\]
\begin{ex}
Let 
$T=
\begin{ytableau}
	1&2&\lg{5}\\
	3&4\\
	6&7\\
\end{ytableau} $ and $B=(1, 3)$.
Then 
\[
\Cu(T, B)=
\begin{ytableau}
	1&2\\
	4&5\\
	6&7\\
\end{ytableau}, \q
\cu(T, B)=3.
\]
\end{ex}

%
%
%
%

Now, let $S\in\ssot_{n}(\l)$.
We will deal with it 
interchangeably 
as a multi-set-valued tableau 
with $\prs=(u_{j})$,
$\tbox\, S= (B_{j})$ 
or $S=(S_{1}, \ds, S_{n})$, 
a sequence of partitions $S_{j}$ with 
$n$ substeps.


\begin{defn}
Let $m\in [n]$. Say 
$m$ is a \kyo{deletion} of $S$ 
if $S_{m-1}\rhd S_{m}$.
Otherwise, $m$ is an \kyo{addition}. 
Let $\del(S)$ and $\add(S)$ denote the set of all deletions and additions of $S$, respectively.
\end{defn}

Observe that
 $m\in\add(S)$ if and only if $B_{m}$ is at an inside corner of $S_{m}$.
 Equivalently, 
 $m\in \del(S)$ if and only if $B_{m}$ is at an outside corner of $S_{m}$.
 
We are now going to define (modified) \kyo{Sundaram  correspondence}
\[
\begin{array}{ccc}
	\th:\ssot_{n}(\l)&\too    &\B(n-|\l|)\ti \ssyt_{}(\l)   \\
	S&   \mto&   (L(S), T(S)).
\end{array}
\]
Let us inductively construct 
a sequence of pairs 
\[
(L_{m}(S), T_{m}(S))\in \B \ti \ssyt{(S_{m})}
\]
for $m=1, \ds, n$ where $\B$ is the set of all 
Burge arrays.

Start with $(L_{1}(S), T_{1}(S))=
\k{\ku, \fb{$u_{1}$}}$. 
Suppose $m\ge2$.
\begin{itemize}
	\item If $m$ is an addition of $S$,
	let 
	$L_{m}(S):=L_{m-1}(S)$
and	
$T_{m}(S)$ the SSYT obtained from $T_{m-1}(S)$  adding the new box $B_{m}$ with entry $u_{m}$. 
This is of shape $S_{m}$.
\item 
Suppose $m$ is a deletion. 
Then $T_{m-1}(S)$ of shape $S_{m-1}$ contains the entry at $B_{m}$, say $x_{m}$. Note that 
$x_{m}$ appeared strictly before the step $u_{m}$ 
(i.e., $x_{m}<u_{m}$). 
Column-unbump $x_{m}$ at $B_{m}$ 
and set $T_{m}(S):=\Cu(T_{m-1}(S), B_{m})$ 
and $v_{m}:=\cu(T_{m-1}(S), B_{m})$ $(\le x_{m}<u_{m}).$
Now let $L_{m}(S)$ be the Burge 
array obtained by rearranging pairs 
in $L_{m-1}(S)$ and \pa{u_{m}}{v_{m}}.
\end{itemize}
Finally, let 
$L(S)=L_{n}(S)$, 
$T(S)=T_{n}(S)$ and define 
\[
\th(S)=(L(S), T(S))\in \B(n-|\l|)\ti \ssyt(\l).
\]
Call $(L(S), T(S))$ the \kyo{Sundaram pair} of 
$S$, $L(S)$ \kyo{Burge part} and $T(S)$ 
\kyo{tableau part}. 

The \kyo{length} of such a pair $(L, T)$ is 
the total number of letters in 
$L$ and $T$. This is nothing but 
the length of $S$ if $(L, T)=(L(S), T(S))$. 
By $\|L\|$ and $\|T\|$, we mean the multisets of 
all entries in such an array and a tableau. 

\begin{lem}[Column Bumping Lemma 
{\cite[p.187]{fu}}]\label{l39}
Let $T$ be an SSYT 
and $v, v' \innn$.
Suppose 
$v'\to (v\to T)$ adds a box $B$ first and $B'$ next. Then
\[
B<_{nE}B'\iff v'\le v.
\]
\end{lem}
\begin{center}

\ytableausetup{mathmode, boxsize=16pt}
\begin{ytableau}
\lg{}&\lg{}&\lg{}&\lg{}&\lg{}&\lg{}\\
\lg{}&\lg{}&\lg{}&\lg{}&{B'}\\
\lg{}&\lg{}&{B}
\end{ytableau}
\end{center}
\begin{lem}
Let 
$\del(S)=\{m_{1}, \ds, m_{r}\}$ 
such that $m_{1}< \cd < m_{r}$ 
and 
$\left(
\pa{u_{m_{l}}}{v_{m_{l}}}\mid 1\le l\le r 
\ri)$ be pairs construted as above. Write $j_{l}=u_{m_{l}}$ and $i_{l}=v_{m_{l}}$ for convenience.
Then, the array 
$L(S)=
\left(\begin{array}{ccccc}
	j_{1}&\cd&j_{r}\\
	i_{{1}}&\cd &i_{{r}}
\end{array}\right)$ 
is Burge.
\end{lem}
\begin{proof}
By construction, 
$j_{1}\le \cd \le j_{r}$, 
$j_{l}=u_{m_{l}}>v_{m_{l}}=i_{{l}}$. 
It is then enough to show that 
$j_{l}=j_{l+1}$ implies 
$i_{l}\le i_{l+1}$.
Suppose 
$j_{l}=j_{l+1}$ 
with 
$j_{l}=u_{m_{l}},
j_{l+1}=u_{m_{l+1}}$, 
$m_{l}, m_{l+1}\in \del(S)$. 
Since a deletion occurs only at an initial segment in a step $u_{m_{l}}$, 
we may assume that 
$m_{l+1}=m_{l}+1$. 
Write $m=m_{l}$ for brevity.
Define 
boxes $B_{m}, B_{m+1}$ and 
letters $v_{m}, v_{m+1}$ by 
\[
T_{m}(S)=\Cu(T_{m-1}(S), B_{m}), 
v_{m}=\cu(T_{m-1}(S), B_{m}),
\]
\[
T_{m+1}(S)=\Cu(T_{m}(S), B_{m+1}), 
v_{m+1}=\cu(T_{m}(S), B_{m+1}).
\]
In other words, 
\[
T_{m-1}(S)=v_{m}\to \k{v_{m+1}\to T_{m+1}(S)}
\]
adds new boxes $B_{m+1}$ first, 
and $B_{m}$ next. 
Notice that 
$B_{m+1}<_{nE} B_{m}$
because these are deletions at consecutive substeps in the same step in $S$. 
By Lemma \ref{l39}, $v_{m}\le v_{m+1}$. 
This means 
$i_{l}\le i_{l+1}$. 
\end{proof}

\begin{ex}
Here we observe 
a semistandard variant of 
\cite[Example 8.8]{sun1}.
\ytableausetup{mathmode, boxsize=25pt}
\[
\te{Let $S=
\begin{ytableau}
	1&244\\
	2&57\\
	346\\
\end{ytableau}\in \ssot_{10}((2, 1, 1))$}
\qte{with 
$\del(S)=\{5, 6, 10\}$.
}
\]

{\renewcommand{\arraystretch}{1.25}
\begin{table}
\caption{$(T_{m}(S))_{m=1}^{10}$}
\label{tab2}
\begin{center}
\fit
{\begin{tabular}{ccccccccccccccc}
	$m$&1&2&3&4&$\mathbf{5}$&$\mathbf{6}$&7&8&9&$\mathbf{10}$\\\h
	$u_{m}$&1&2&2&3&4&4&4&5&6&7\\\h
	$T_{m}(S)$&$\begin{array}{ccccc}
		1\\
		\mb{}\\\mb{}
	\end{array}$&
$\begin{array}{ccccc}
		1\\
		2\\\mb{}
	\end{array}$&
	$\begin{array}{ccccc}
		1&2\\
		2\\\mb{}
	\end{array}$&
	$\begin{array}{ccccc}
		1&2\\
		2\\3\\\mb{}
	\end{array}$&
	$\begin{array}{ccccc}
		1\\
		2\\3\\\mb{}
	\end{array}$&
	$\begin{array}{ccccc}
		1\\
		2\\\mb{}
	\end{array}$&
	$\begin{array}{ccccc}
		1&4\\
		2\\\mb{}
	\end{array}$&
	$\begin{array}{ccccc}
		1&4\\
		2&5\\\mb{}
	\end{array}$
	&
	$\begin{array}{ccccc}
		1&4\\
		2&5\\6
	\end{array}$&
		$\begin{array}{ccccc}
		1&4\\
		5\\6\\\mb{}
	\end{array}$
\end{tabular}}
\end{center}
\end{table}}
Table \ref{tab2} shows $(T_{m}(S))_{m=1}^{10}$.

Observe that 
\[
L_{5}(S)=
\left(\begin{array}{ccccc}
	4\\
	2
\end{array}\right), 
L_{6}(S)=
\left(\begin{array}{ccccc}
	4&4\\
	2&3
\end{array}\right), 
L_{10}(S)=
\left(\begin{array}{ccccc}
	4&4&7\\
	2&3&2
\end{array}\right)
\]
and hence 
\ytableausetup{mathmode, boxsize=15pt}
\[
\th(S)=
\k{
\left(\begin{array}{ccccc}
	4&4&7\\
	2&3&2
\end{array}\right), 
\begin{ytableau}
	1&4\\
	5\\
	6\\
\end{ytableau}
}.
\]
\end{ex}


\begin{lem}
$\th:\ssot_{n}(\l)\to \B(n-|\l|)\ti \ssyt_{}(\l)
$ is a weight-preserving bijection.
\end{lem}
\begin{proof}
We prove the claim by induction on $n\ge 1$. 
First, we show that $\th$ is injective.
If $n=1$, then 
$\th\k{\fb{$u_{1}$}}=\k{\ku, \fb{$u_{1}$}}$ is certainly so. Suppose $n\ge2$ and the assertion holds for $n-1$. 
Let 
$S, S'\in \ssot_{n}(\l)$,
\[
\pr S=(u_{j}),\q
\pr S'=(u'_{j}),\q 
\te{box}\, S=(B_{j}) 
\te{ and }
\te{box}\, S'=(B'_{j})
\]
and 
assume that 
$\th(S)=\th(S')$.
Both of $u_{n}$ and $u_{n}'$ 
are maxima of 
\[
\|L(S)\|\cup \|T(S)\|=
\|L(S')\|\cup \|T(S')\|.
\]
Thus, they must be equal. Say $x=u_{n}=u_{n}'$ for convenience. 
Let $\wti{S}$ and $\wti{S'}$ be the 
SSOTs 
obtained by deleting $x$ from $S$ and $S'$ 
expressed as set-valued tableaux, respectively.
\begin{enumerate}
\item Suppose $n\in \add(S)$ so that 
$u_{n}$ appears at $B_{n}$ in $T(S)$. 
This $B_{n}$ must be the rightmost box in $T(S)$  which contains $x$. Since $x$ appears also in $T(S')$ (because $T(S')=T(S)$), there exists 
$m\in \add(S')$ such that $u_{m}'=x$. 
Now, $x=u_{n}'\ge u_{m}'=x$ implies 
$n\in \add(S')$, too. 
Again, 
$u_{n}'$ must be in the rightmost box in $T(S')$ with $x$. Thus, $B_{n}=B_{n}'$. 
By removing $x$ at this box from $T(S)=T(S')$, we have 
\[
\th(\wts)=(L(\wts), T(\wts))=
(L(S), T(\wts))=(L(S')), T(\wti{S'}))=
(L(\wti{S'})), T(\wti{S'}))=
\th(\wti{S'}).
\]
By induction, we have $\wts=\wti{S'}$ and hence $S=S'$.
\item Suppose $n\in \del(S)$. 
Then $u_{n}'$ cannot appear in $T(S')(=T(S))$, 
i.e., 
$n\in \del(S')$.
The rightmost pair 
in $L(S)$ is 
$\ret{u_{n}}{i_{n}}=
\ret{u_{n}'}{i'_{n}}$
and so $i_{n}=i_{n}'$. 
It follows from 
\[
T(\wts)=
\Cu(T(S), B_{n})=
\Cu(T(S'), B_{n}')=T(\widetilde{S'}) 
\te{ and }
\]
\[
L(\wts)=
L(S)\sm \ret{u_{n}}{i_{n}}
=
L(S')\sm \ret{u'_{n}}{i'_{n}}
=L(\wti{S'})
\]
that 
$\th(\wts)=\th(\wti{S'})$. By induction, 
we have 
$\wts=\wti{S'}$ and hence $S=S'$.
\end{enumerate}

Next, we show that 
$\th$ is surjective also by induction on $n$.
For $n=1$, it is surjective as seen above.
Suppose $n\ge2$ and the assertion holds for $n-1$. 
Let 
$(L, T)\in \B(n-|\l|)\ti \ssyt(\l)$ and 
it is enough to find $S$ such that $\th(S)=(L, T)$.
Let $x=\max (\|L\|\cup \|T\|)$.
\begin{enumerate}
	\item Suppose $x\in \|T\|$. 
	Let $B_{n}$ be the rightmost box in $T$ which contains $x$. 
	Set $\wti{L}=L$, 
	$\wti{T}=\Cu(T, B_{n})$, 
	say $\wti{\l}=\sh \wti{T}$. 
	For the pair 
	\[
(\wti{L}, \wti{T})
	\in \B(n-1-|\wti{\l}|)\ti \ssyt(\wti{\l}),\]
the inductive hypothesis claims that 
there exists $\wti{S}\in \ssot_{n-1}(\wti{\l})$ such that $\th(\wts)=(\wti{L}, \wti{T})$.
Let $S$ be the SSOT with all entries of 
$\wti{S}$ and $u_{n}=x$ at $B_{n}$. 
We have 
$L(S)=\wti{L}=L$ 
and $T(S)=\wti{T}\oplus(B_{n}, x)=T$ so that 
\[
\th(S)=(L(S), T(S))=
(L, T).
\]
\item Suppose $x\in \|L\|\setminus\|T\|$. 
Say the rightmost pair of $L$ is 
$\ret{x}{i_{n}}$ and $u_{n}:=x$.
Let $\wti{L}$ be the array obtained from $L$ by removing one $\ret{x}{i_{n}}$.
Let 
$\wti{T}=i_{n}\to T$. 
Say the new box is $B_{n}$ and 
$\sh \wti{T}=\wti{\l}$. By induction, 
there exists $\wti{S}\in \ssot_{n-1}(\wti{\l})$ such 
that 
$\th(\wts)=(\wti{L}, \wti{T})$. 
Let $S$ be the SSOT with entries of 
$\wti{S}$ and $u_{n}=x$ at $B_{n}$.
By construction, 
\[
L(S)=\wti{L}\oplus \ret{x}{i_{n}}=L
\qte{and}\q
T(S)=\Cu(\wti{T}, B_{n})=T
\]
which implies $\th(S)=(L, T)$.
\end{enumerate}
Finally, it is clear that 
$\th$ preserves the weight 
because $\|S\|=\|L(S)\|\sqcup \|T(S)\|$. 
\end{proof}

\begin{thm}\label{t31}
We have 
\[
\ss_{\l}(\x)=
\prod_{i<j}(1-x_{i}x_{j})^{-1}s_{\l}(\x).
\]
\end{thm}
\begin{proof}
The bijectivity of $\th$ shows 
\[
\ss_{\l}(\x)=
\k{\sum_{\bev}s_{\b}(\x)}
s_{\l}(\x).
\]
Our claim follows from this and 
the Littlewood identity
\[
\prod_{i<j}(1-x_{i}x_{j})^{-1}=
\sum_{\bev}s_{\b}(\x).
\]
\end{proof}
\begin{cor}
$\ss_{\l}(\x)$ is symmetric.
\end{cor}
\begin{rmk}
In \cite[Theorem 4.28]{km}, the first and second authors 
proved such symmetry 
under a finite-variable setting. 
Here we include the case with infinite ones.
\end{rmk}
\np

\np

\section{Applications}

\subsection{Schur positivity}

Next, we will prove the Schur positivity, the orthogonality and the 
linear independence of SSOT functions in analogy of the Schur functions. 
Note that the degree of $\ss_{\l}(\x)$ is \kyo{unbounded}. 
It is thus technically easier to deal with its  homogeneous part $\ss_{\l, n}(\x)$. 
We will also see an application to 
\kyo{polytopes} which play an important role in  combinatorics.
Toward Theorems \ref{t41}, \ref{t42} and \ref{t43}, we introduce some vocabularies and notations.





\begin{defn}
A word $w=w_{1}\cd w_{n}$ of positive integers 
 is \eh{reverse Yamanouchi} if 
the weight of the subword $w_{m}w_{m+1}\cd w_{n}$ for every $m\le n$ is a partition.
Let $\l, \m$ and $\n$ be partitions.
A tableau $S$ of skew shape $\n/\l$ 
is a \kyo{Littlewood-Richardson (LR) tableau} of weight $\m$ if the row word of $S$ is reverse Yamanouchi of weight $\m$.
Let $\te{LR}_{\l\m}^{\n}$ 
denote 
the set of such tableaux and 
$\clmn=
|\te{LR}_{\l\m}^{\n}|.$
\end{defn}
Note that 
$\clmn$ is 0 unless 
$\l, \m\sub \n$ and $|\n|=|\l|+|\m|$.

\begin{fact}[LR rule]
We have 
\[
s_{\l}s_{\mu}=\sum_{\nu}c_{\l\m}^{\n} s_{\n}.
\]
In particular, $c_{\l\m}^{\n}=c_{\m\l}^{\n}$.
\end{fact}

A skew shape is a 
\kyo{vertical strip} if 
it contains at most one box in each row.

\begin{defn}
For $n\in N(\l)$, define $V^{n}(\l)$ to be 
the set of all partitions $\n$ which satisfy the following.  
\begin{enumerate}
	\item $|\n|=n$.
	\item $\n$ is obtained 
from $\l$ by adding several 
vertical strips of even size. 
\end{enumerate}
Set $V(\l)=\bigsqcup_{n\in N(\l)}V^{n}(\l)$.
\end{defn}
Observe that 
$V^{n}(\l)\ne \ku$ whenever 
$n\in N(\l)$
and in particular, 
$V^{|\l|}(\l)=\{\l\}$. 
The number of vertical strips of even size 
that can be added to $\l$ 
to obtain $\nu\in V^{n}(\l)$ is at most 
$(n-|\l|)/2$.


\begin{ex}
All shapes in $V^{7}\k{
\ytableausetup{mathmode, boxsize=8pt}\yd{2,1}}$ are:
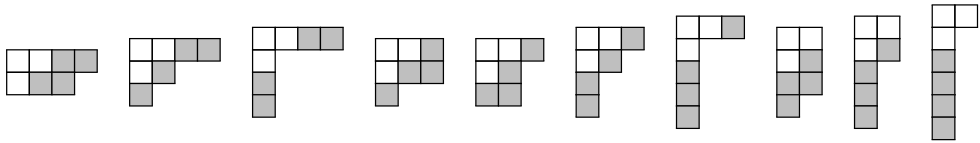
\[
\begin{ytableau}
	\mb{}&\mb{}&\lg&\lg\\
	\mb{}&\lg&\lg\\
\end{ytableau}\q
\begin{ytableau}
	\mb{}&\mb{}&\lg&\lg\\
	\mb{}&\lg\\
	\lg\\
\end{ytableau}\q
\begin{ytableau}
	\mb{}&\mb{}&\lg&\lg\\
	\mb{}\\
	\lg\\
	\lg\\
\end{ytableau}\q 
\begin{ytableau}
	\mb{}&\mb{}&\lg\\
	\mb{}&\lg&\lg\\
	\lg\\
\end{ytableau}\q
\begin{ytableau}
	\mb{}&\mb{}&\lg\\
	\mb{}&\lg\\
	\lg&\lg\\
\end{ytableau}
\q
\begin{ytableau}
	\mb{}&\mb{}&\lg\\
	\mb{}&\lg\\
	\lg\\
	\lg\\
\end{ytableau}\q
\begin{ytableau}
	\mb{}&\mb{}&\lg\\
	\mb{}\\
	\lg\\\lg\\\lg\\
\end{ytableau}\q
\begin{ytableau}
	\mb{}&\mb{}\\
	\mb{}&\lg\\
	\lg&\lg\\
	\lg\\
\end{ytableau}
\q
\begin{ytableau}
	\mb{}&\mb{}\\
	\mb{}&\lg\\
	\lg\\
	\lg\\
\lg
\end{ytableau}\q
\begin{ytableau}
	\mb{}&\mb{}\\
	\mb{}\\
	\lg\\
	\lg\\
	\lg\\
\lg
\end{ytableau}
\]
\end{ex}

We will prove the following theorem below.
\begin{thm}
\label{t41}
Let $n\in N(\l)$. 
Then we have 
\[
\ss_{\l, n}=
\sum_{
\substack{\b' \te{ even}
\\
\n\in V^{n}(\l)}
}
c_{\b\l}^{\n}s_{\n}.
\]
In particular, this is a finite sum of Schur functions 
and hence Schur-positive.
\end{thm}


\begin{lem}\hf
\label{l46}
\begin{enumerate}
\item 
	Let $U(\l)$ be the tableau of shape $\l$ such that its row $i$ consists of only $i$. 
Then 
\[
c_{\m\l}^{\n}=
|\{T\in \ssyt(\m)|T\cdot U(\l)=U(\n)\}|
\]
where $\cdot$ means the product of 
SSYTs as in \cite[Chapters 1, 2]{fu}.
\item 
Let 
$v=v_{m}\cd v_{1}$ be a column word ($v_{i+1}>v_{i}$) and $T$ any SSYT of shape $\l$. 
Then column insertions 
$v_{m}\to (\cd \to (v_{1}\to T))$ 
adds a vertical strip to $\l$.
\end{enumerate}
\end{lem}
\begin{proof}
See \cite[p.66 and p.187, Column Bumping Lemma]{fu}.
\end{proof}

\begin{lem}\label{l47}
If $c_{\b\l}^{\n}\ne 0$ for some 
$\b$ such that $\b'$ is even 
and $n=|\n|$, 
then $\n\in V^{n}(\l)$.
\end{lem}
\begin{proof}
Suppose 
$c_{\b\l}^{\n}\ne 0$ with $\b, \n$ as above.
Thanks to Lemma \label{l46} (1), there exists some tableau $T$ of shape $\b$ such that $T\cdot U(\l)=U(\n)$.
Lemma \label{l46} (2) says that column-inserting a column word 
of each column of $T$ (of even length) 
adds a vertical strip of even size 
to $\l$. 
This implies $\n\in V^{n}(\l)$.
\end{proof}

\begin{proof}[Proof of Theorem \ref{t41}]
We have 
\[
\ss_{\l, n}=
\k{\prod_{i<j}(1-x_{i}x_{j})^{-1}}
s_{\l}
=\sum_{\b' \te{ even}} s_{\b}s_{\l}
=\sum_{\b' \te{ even}}
\sum_{\n} c_{\b\l}^{\n}s_{\n}
=
\sum_{
\substack{\b' \te{ even}
\\
\n\in V^{n}(\l)}
}
c_{\b\l}^{\n}s_{\n}
\]
where the last equality follows from 
Lemma \label{l47}.
\end{proof}

\ytableausetup{smalltableaux}
\begin{ex}
Suppose that $\l=
\yd{2, 1}$, $\n\in V^{5}(\l)$ 
and $c_{\b\l}^{\n}\ne 0$, $\b'$ even. 
In such a case, $|\b|$ must be $2$. 
Thus, the only choice of $\b$ is $(1, 1)$
and $T$ of shape $\b$ such that 
$T\cdot U(\l)=U(\n)$ 
is 
\ytableausetup{mathmode, boxsize=13pt}
$\begin{ytableau}
	1\\2
\end{ytableau}$ 
(the row word of $T$ must be a reverse Yamanouchi  word of weight (1, 1))
so that $c_{\b\l}^{\n}=1$.
We thus have 
\ytableausetup{mathmode, boxsize=8pt}
\[
\ss_{\,
\begin{ytableau}
	\mb{}&\mb{}\\
	\mb{}
\end{ytableau}, 5
}
=
s_{\,
\begin{ytableau}
	\mb{}&\mb{}&\lg{}\mb{}\\
	\mb{}&\lg{}\mb{}
\end{ytableau}
}
+
s_{\,
\begin{ytableau}
	\mb{}&\mb{}&\lg{}\mb{}\\
	\mb{}\\
	\lg{}\mb{}
\end{ytableau}
}+
s_{\,
\begin{ytableau}
	\mb{}&\mb{}\\
	\mb{}&\lg{}\mb{}\\
	\lg{}\mb{}
\end{ytableau}
}
+
s_{\,
\begin{ytableau}
	\mb{}&\mb{}\\
	\mb{}\\
	\lg{}\mb{}\\
	\lg{}\mb{}\\
\end{ytableau}
}.
\]
\end{ex}

Let 
$\L^{n}=\k{\cc[[\x]]^{n}}^{S_{\mg}}$ 
be the space of homogeneous 
symmetric functions of degree $n$. 
For a nonnegative integer $m$, 
let 
\[
N(m)=\{n\mid n\ge m, \hou{n}{m}{2}\}.
\]
Let $\parm$ denote the set of all partitions of 
$m$.
There exists 
the inner product, \kyo{Hall inner product}, 
$\la\ph{a},\ph{b} \ra$ on 
$\L^{n}$ such that 
$(s_{\l}|\l \in \te{Par}(n))$ forms 
an orthonormal basis:
\[
\la s_{\l}, s_{\m}\ra=\delta_{\l\m} 
\qte{\cite[p.78]{fu}}.
\]
Here we would like to know if we can say something similar on $\ss_{\l, n}$. 
As below, we prove results on partial orthogonality and linear independence.
For $\l$ and $\m\in \parm$, define 
\[
V^{n}(\l, \m)=
V^{n}(\l)\cap V^{n}(\m).
\]
Let us say that $\l$ and $\m$ are $n$-similar 
if $V^{n}(\l, \m)\ne \ku$.
Otherwise, we say 
they are $n$-orthogonal.

\begin{lem}\label{l48}
Let $\l, \m\in \parm$ and 
$n_{0}, n\in N(m)$ with $n_{0}<n$.
If $\l$ and $\m$ are $n_{0}$-similar, 
then they are also $n$-similar.
\end{lem}
\begin{proof}
Suppose 
$\n_{0}\in V^{n_{0}}(\l, \m)$.
To $\n_{0}$, add several 
vertical strips of even size to construct 
some partition $\n\in \te{Par}(n)$. 
This immediately implies $\n\in V^{n}(\l, \m)$.
\end{proof}

\begin{lem}\label{l49}
For any $\l$ and $\m\in \parm$, 
we have 
$V^{n}(\l, \m)\ne \ku$ for some $n\in N(m)$.
\end{lem}
\begin{proof}
Certainly, $m^{2}\in N(m)$ 
since $m^{2}-m$ is even.
Adding horizontal strips of even size to $\l, \m$ respectively, we can construct $(m^{m})$.
\end{proof}

\begin{lem}\label{l410}
Let $\l, \m\in \parm$ and $n\in N(m)$. Then we have 
\[
\la\ss_{\l, n}, \ss_{\m, n}\ra=
\sum_{
\substack{\a',\,\b'\te{ even}\\\n}
}
{c_{\a\l}^{\n}c_{\b\m}^{\n}}\,\,(\ge 0).
\]
Moreover, this integer is positive $\iff$ 
$V^{n}(\l, \m)\ne \ku$.
\end{lem}
\begin{proof}
Observe that 
\[
\la\ss_{\l, n}, \ss_{\m, n}\ra=
\left\la
\sum_{
\substack{\a'\te{ even}\\\n_{1}}
}
c_{\al \l}^{\n_{1}}s_{\n_{1}}, 
\sum_{
\substack{\b'\te{ even}\\\n_{2}}
}
c_{\b \m}^{\n_{2}}s_{\n_{2}}
\ri\ra.
\]
Expand this by the linearity of 
$\la\ph{a}, \ph{a}\ra$.
Then it follows from 
the 
orthonormality 
of Schurs that we have 
\[
\la\ss_{\l, n}, \ss_{\m, n}\ra=
\sum_{
\substack{\a',\,\b'\te{ even}\\\n}
}
{c_{\a\l}^{\n}c_{\b\m}^{\n}}.
\]

Moreover, 
$c_{\a\l}^{\n}c_{\b\m}^{\n}>0$ holds 
for some partitions $\a, \b$ and $\n$ 
such that both of $\a'$ and $\b'$ are even 
if and only if $\n \in V^{n}(\l, \m).$
\end{proof}

\begin{thm}\label{t42}
Let $\l, \m\in \parm$.
There exists a unique $n_{0}\in N(m)$ 
such that 
\begin{align*}
	\la \ss_{\l, n}, \ss_{\m, n}\ra&=0\qte{for all 
	$n\in N(m)$ such that $n< n_{0}$},
	\\\la \ss_{\l, n}, \ss_{\m, n}\ra&
	>0\qte{for all $n\in N(m)$ such that 
	$n\ge n_{0}$}.
\end{align*}
\end{thm}
\begin{proof}
By Lemmas \ref{l48}, \ref{l49} and \ref{l410}, we can choose $n_{0}=\min\{n\in N(m)\mid V^{n}(\l, \m)\ne \ku\}$.
\end{proof}

Note that $\la s_{\l}, s_{\m}\ra=\delta_{\l\m}$ is merely a special case of this.
\begin{ex}\hf
Let $\l=
\yd{3}$ and $\mu=\yd{1,1,1}$. 
Observe that 
$\l$ and $\m$ are 5-orthogonal 
because 
\[
V^{5}(\l, \m)=
\left\{
\begin{ytableau}
	\mb{}&\mb{}&\mb{}&\lg{}\\
	\lg{}\mb{}\\
\end{ytableau},
\begin{ytableau}
	\mb{}&\mb{}&\mb{}\\
	\lg{}\mb{}\\
	\lg{}
\end{ytableau}
\ri\}\cap
\left\{
\begin{ytableau}
	\mb{}&\lg{}\\
	\mb{}&\lg{}\\
	\mb{}\\
\end{ytableau}, 
\begin{ytableau}
	\mb{}&\lg{}\\
	\mb{}\\
	\mb{}\\
	\lg{}\\
\end{ytableau},
\begin{ytableau}
	\mb{}\\
	\mb{}\\
	\mb{}\\
	\lg{}\\
	\lg{}
\end{ytableau}
\ri\}
=\ku.
\]
However, 
they are 7-similar because 
\[
V^{7}(\l, \m)=
\left\{
\yd{3,2,2}, \yd{3,2,1,1},
\yd{3,1,1,1,1}
\ri\}
\ne\ku.
\]
\end{ex}



We spend the rest of this section to prove the following theorem.

\begin{thm}\label{t43}
$(\ss_{\l, n}|\l \in \te{Par}(m))$ is linearly 
independent.
\end{thm}

Our proof follows from a typical argument using 
the \kyo{dominance order}. 
For $\l, \m\in \parm$, 
define $\l\ge \mu$ if 
\[
\l_{1}+\cd+\l_{i}\ge \m_{1}+\cd+\m_{i}
\]
for all $i$ (adding 0's to the end if necessary). 

For a partition 
$\l=(\lam_{1}, \ds, \lam_{l})$ with $\l_{i}\ge 0$,  $n\in N(\l)$, 
say $n=|\l|+2r$, define 
$\ol{\l}^{\,n}\in V^{n}(\l)$ by 
\[
\ol{\l}^{\,n}=(\l_{1}+r, \l_{2}+r, \l_{3},\ds, \l_{l}).
\]
Below, we drop $n$ from notation for simplicity 
and write $\ol{\l}=\ol{\l}^{\,n}$.

\begin{lem}\label{l415}
Let 
$\l\in \parm$ and $n\in N(\l)$.
We have 
$\max V^{n}(\l)=\oll$ where $\max$ is with respect to the dominance order. 
\end{lem}

\begin{proof}
Suppose $\n\in V^{n}(\l)$.
Say 
\[
\l=(\l_{1}, \ds, \l_{l}), 
\n=(\n_{1}, \ds, \n_{l'}), l'\ge l, n=|\l|+2r.
\]
By the 
definition of $V^{n}(\l)$, 
there exist $r_{1}, \ds, r_{l'}\ge0$ such that 
\[
\n_{i}=\l_{i}+r_{i}, r_{i}\le r \qte{for all $i$}
\]
and $r_{1}+\cd+r_{l'}=2r$. 
It follows that 
\[
\oll_{1}=
\l_{1}+r\ge \l_{1}+r_{1}=\n_{1},
\]
\[
\oll_{1}+\oll_{2}=
\l_{1}+r+ \l_{2}+r\ge
\n_{1}+r_{1}+ \n_{2}+r_{2}
=\n_{1}+\n_{2},
\]
and 
for $i\ge 3$,
\[
\oll_{1}+\cd+\oll_{i}
= \l_{1}+\cd+\l_{i}+2r
\ge \l_{1}+\cd+\l_{i}+r_{1}+\cd+r_{i}
=\n_{1}+\cd+\n_{i}.
\]
This proves that 
$\oll$ is the maximum element in $V^{n}(\l)$ with  respect to the dominance order. 
\end{proof}

\begin{ob}
Let $\l, \m\in \parm$ and $n\in N(m)$.
\begin{enumerate}
	\item $\l=\m $ if and only if $\oll=\olm.$
	\item $\l>\m$ if and only if 
	 $\oll>\olm$.
\end{enumerate}
Thus, 
$\l\lr \oll$ gives an order-preserving bijection 
between $\parm$ and 
\[
\ol{\parm}:=
\{\ol{\l}\in \parn\mid \l\in \parm\}.
\]
\end{ob}


\begin{lem}\label{l416}
${\sum_{\bev}c_{\b\l}^{\oll}}=1.$
In particular, if $n\in N(\l)$, then 
we have 
\[
\ss_{\l, n}=
s_{\oll}+
\sum_{
\substack{\bev\\\n<\oll
}
}c_{\b\l}^{\n}s_{\n}.
\]
\end{lem}
\begin{proof}
Let $\b$ be 
a partition such that 
$\b'$ is even and $c_{\b\l}^{\oll}\ne 0$.
There exists a tableau $T$ of shape $\b$ such that  
$T\cdot U(\l)=U(\oll)$ (Lemma \ref{l46}).
It follows that 
\[
\wt\, T=\wt \oll-\wt \l=(r, r, 0, \ds, 0).
\]
Since $\b'$ is even, the only possibility 
is 
$T=
\ytableausetup{mathmode, boxsize=16pt}
\begin{ytableau}
	1&\cd&1\\
	2&\cd&2\\
\end{ytableau}$ 
with $\b=(r, r)$, too. 
The second assertion holds due to 
Lemma \ref{l415}.
\end{proof}

\begin{proof}[Proof of Theorem \ref{t43}]
Let $M=|\parm|$
and choose distinct 
$\l^{(1)}, \ds, 
\l^{(M)}\in \parm$ such that 
$\l^{(i)}$ is maximal 
in the subposet 
$\{\l^{(i)}, \ds, \l^{{(M)}}\}$ 
with 
$\l^{(1)}=(m)=\max \parm$.
To show the 
linear independency, 
suppose 
\[
\dsum_{i=1}^{M}b_{i}\ss_{\ol{\l^{(i)}}, n}(\x)=0	, \q b_{i}\in \cc.
\]
Expand this by Schur functions.
By Lemma \ref{l416}, the coefficient of 
$s_{\,\ol{\l^{(1)}}}$ comes only 
from $\ss_{\,\ol{\l^{(1)}}, n}$, 
that is, 
\[
b_{1}
s_{\,\ol{\l^{(1)}}}+
\te{(all other terms without $s_{\,\ol{\l^{(1)}}}$)}=0.
\]
Since 
$(s_{\n} |\n\in \parn)$ forms a basis of $\Lambda^{n}$, $b_{1}$ must be $0$.
Inductively, we can show $b_{2}=\cd=b_{M}=0$.
\end{proof}

\subsection{Newton polytope}

Let 
$f(x_{1}, \ds, x_{k})=\sum_{}c_{a}x^{a}
\in \zz[x_{1}, \ds, x_{k}]$ be a polynomial. Its support is 
\[
\supp f=\{a\in \zz_{\ge0}^{k}\mid 
c_{a}\ne 0\}.
\]
The \kyo{Newton polytope} of $f$ is 
\[
\te{New}(f)=
\te{conv}(\supp f)
\,\k{\sub \mathbf{R}^{k}}
\]
where \te{conv} means the convex hull.
\begin{defn}
Say $f$ has \kyo{Saturated Newton polytope} (SNP) 
if 
\[
\te{New}(f)\cap \zz^{k}_{\ge0}=\supp f.
\]
\end{defn}

For example, all 
Schur polynomials have SNP while 
not all Schur-positive polynomials do \cite{mty}.
It is then natural to study the SNP property 
of SSOT polynomials as well.
\begin{thm}\label{t44}
If $n\in N(\l)$, then $\ss_{\l, n}(x_{1}, \ds, x_{k})$  has SNP.
\end{thm}
For a proof, our discussions in the previous subsection play a role with this lemma.

\begin{lem}[{\cite[Proposition 2.5]{mty}}]\label{l418}
Suppose that a polynomial 
$f(x_{1}, \ds, x_{k})$ satisfies 
\[
f=\sum_{\n\in \parn} b_{\n}s_{\n}
\]
such that for some $\oll\in \parn$, 
$b_{\oll}\ne 0$ and moreover 
$b_{\n}\ne 0$ only if 
$\n\le \oll$. 
Then $f$ has SNP.
\end{lem}

\begin{proof}[Proof of Theorem \ref{t44}]
By Lemma \ref{l416}, we find that 
\[
f=\ss_{\l, n}=
{\sum_{
\substack
{\beta' \text{even}\\
\n\in V^{n}(\l), \n\le \oll}
}
\underbrace{c_{\b\l}^{\n}}_{b_{\n}}}
s_{\n}
\]
satisfies all assumptions of Lemma 
\ref{l418}.
\end{proof}

\end{document}